\documentclass{amsart}
\usepackage{amsmath,amssymb,latexsym,graphics,verbatim}
\usepackage[margin=2cm]{geometry}
\theoremstyle{plain}

\theoremstyle{definition}

\theoremstyle{remark}

\begin{document}
\title{A Proof of the Thue-Siegel Theorem
\\about the Approximation of Algebraic Numbers for Binomial Equations}
\author{Kurt Mahler, translated by Karl Levy}
\date{\today}

\maketitle
\pagestyle{myheadings}\markboth{KURT MAHLER, Translated by Karl Levy}{A PROOF OF THE THUE-SIEGEL THEOREM {\tiny about the Approximation of Algebraic Numbers for Binomial Equations}}

In 1908 Thue (1) showed that algebraic numbers of the special form $\xi =\sqrt[n]{\frac{a}{b}}$ can, for every positive $\epsilon$, only be sharply approximated by finitely many rational numbers $\frac{p}{q}$ with the following inequality holding
\[
\left|\xi -\frac{p}{q}\right|\leq q^{-(\frac{n}{2}+1+\epsilon)}
.\]
The proof uses, if perhaps in a somewhat hidden way, the continued fraction expansion of the binomial series $(1-z)^\omega$.  In further work about the approximation of algebraic numbers (2,3) famously Thue used instead a completely different tool, the drawer method  of Dirichlet, and showed further that the above statement holds for any algebraic number. Thue's methods were later generalized by Siegel(4,5,6,7) who showed, among other things, that for every algebraic number in the above inequality the exponent $\frac{n}{2}+1+\epsilon$ could be replaced by $\frac{n}{m}+m-1+\epsilon$, where $m$ is some natural number.

This note demonstrates a generalization of Thue's methods in (1);  Like Thue I restrict myself to the roots $\xi =\sqrt[n]{\frac{a}{b}}$ of the binomial equations.  The continued fraction expansion of the binomial series is generalized and algebraic approximation functions are given instead of rational approximation functions.  In doing so I proceed exactly as in my work on the exponential function(8).  Integrals are set up for the approximation functions; thus the estimates become much easier and you can  prove Thue's theorem with Siegel's Exponents for the binomial algebraic equations without difficulty and without use of the pigeonhole principle.
\section*{I.}
\subsection*{1}
Let $\varrho_1, \varrho_2,\dots,\varrho_m$ be $m$ natural numbers and let $\omega_1,\omega_2,\dots, \omega_m$ be $m$ complex numbers such that no pairwise difference 
\[\omega_h-\omega_k \qquad (h,k = 1,2,\dots ,m);h\neq k\]
is an integer.  From known theorems about homogeneous linear equations there are $m$ polynomials
\[
A_k \left( z \,\,\vline
\begin{tabular}{ccc}
$\omega_1$ & $\dots$ & $\omega_m$ \\
$\varrho_1$ & $\dots$ & $\varrho_m$
\end{tabular} 
\right)
\qquad
(k=1,2,\dots ,m)
\]
that do not simultaneously and identically vanish and that are respectively of degree at most
\[\varrho_1 -1, \varrho_2 -1,\ldots,\varrho_m -1,\]
 so that in the power series expansion of the expression
\[
\sum_{k=1}^m
A_k \left( z \,\,\vline
\begin{tabular}{ccc}
$\omega_1$ & $\dots$ & $\omega_m$ \\
$\varrho_1$ & $\dots$ & $\varrho_m$
\end{tabular} 
\right)
(1-z)^{\omega_k}
=
\sum_{l=0}^\infty a_l z^l
\]
all coefficients $a_l$ with
\[0\leq l<\varrho_1 + \varrho_2 +\dots+ \varrho_m -1\]
are zero.  We rewrite these expressions as
\[ 
R\left( z\,\,\vline
\begin{tabular}{ccc}
$\omega_1$ & $\dots$ & $\omega_m$ \\
$\varrho_1$ & $\dots$ & $\varrho_m$
\end{tabular} \right) 
=
\sum_{k=1}^m A_k\left(z\,\,\vline
\begin{tabular}{ccc}
$\omega_1$ & $\dots$ & $\omega_m$ \\
$\varrho_1$ & $\dots$ & $\varrho_m$
\end{tabular} \right) 
(1-z)^{\omega_k}
.\]
Then it can easily be shown that
\[\frac{d^\varrho_1}{dz^\varrho_1}
\left\{
(1-z)^{-\omega_1}
R\left( z \,\,\vline
\begin{tabular}{ccc}
$\omega_1$ & $\dots$ & $\omega_m$ \\
$\varrho_1$ & $\dots$ & $\varrho_m$
\end{tabular} \right)
\right\}\]
is the same as
\[
R\left( z \,\,\vline
\begin{tabular}{ccc}
$\omega_2-\omega_1-\varrho_1$ & $\dots$ & $\omega_m-\omega_1-\varrho_1$ \\
$\varrho_2$ & $\dots$ & $\varrho_m$
\end{tabular} \right)
\]
and thus that consequently the coefficient of the $(\varrho_1 + \varrho_2 +\dots+ \varrho_m -1)$-th power of $z$ in the power series expansion of $R\left(z\,\,\vline\begin{tabular}{ccc}$\omega_1$ & $\dots$ & $\omega_m$ \\$\varrho_1$ & $\dots$ & $\varrho_m$ \end{tabular} \right)$ is not equal to zero.  We choose the coefficient to be:

\[
\frac{\Gamma(\varrho_1)\ldots\Gamma(\varrho_m)}{\Gamma(\sigma)}
\qquad
\left(\sigma=\sum_{k=1}^m\varrho_k\right),
\]
and thus $R\left(z\,\,\vline\begin{tabular}{ccc} $\omega_1$ & $\dots$ & $\omega_m$ \\ $\varrho_1$ & $\dots$ & $\varrho_m$ \end{tabular} \right)$ is \textit{uniquely} determined.  Below $R\left(z\,\,\vline\begin{tabular}{ccc} $\omega_1$ & $\dots$ & $\omega_m$ \\ $\varrho_1$ & $\dots$ & $\varrho_m$  \end{tabular} \right)$ is always to be understood as such\footnote{See the proof in my paper (8) wherein the analogous considerations for the exponential function are gone through completely.}.

Thus we have the identity 
\begin{align*}
\left\{\frac{d^{\varrho_{m-1}}}{dz^{\varrho_{m-1}}}(1-z)^{\omega_{m-2}+\varrho_{m-2}-\omega_{m-1}} \right\}
\left\{ \frac{d^{\varrho_{m-2}}}{dz^{\varrho_{m-2}}}(1-z)^{\omega_{m-3}+\varrho_{m-3}-\omega_{m-2}} \right\}\cdots \\
 \left\{ \frac{d^{\varrho_{1}}}{dz^{\varrho_{1}}}(1-z)^{-\omega_{1}} \right\}\frac{R\left(z\,\,\vline\begin{tabular}{ccc} $\omega_1$ & $\dots$ & $\omega_m$ \\ $\varrho_1$ & $\dots$ & $\varrho_m$  \end{tabular} \right)}{\Gamma(\varrho_1)\dots\Gamma(\varrho_m))}=\frac{R\left(z \,\,\vline \begin{tabular}{c} $\omega_{m} -\omega_{m-1}- \varrho_{m-1}$ \\ $\varrho_{m} $ \end{tabular} \right)}{\Gamma(\varrho_m)}
\end{align*}
and since we have
\[
\frac{R\left( z \,\,\vline \begin{tabular}{c} $\omega_{m} -\omega_{m-1}- \varrho_{m-1}$ \\ $\varrho_{m} $ \end{tabular} \right)}{\Gamma(\varrho_m)} 
=
 \frac{z^{\varrho_m -1}(1-z)^{\omega_{m}-\omega_{m-1}-\varrho_{m-1}}}{\Gamma(\varrho_m)} 
\]
we can write
\begin{align*}
\frac{R\left(z\,\,\vline\begin{tabular}{ccc} $\omega_1$ & $\dots$ & $\omega_m$ \\ $\varrho_1$ & $\dots$ & $\varrho_m$ \end{tabular}\right)}{\Gamma(\varrho_1)\ldots\Gamma(\varrho_m)}
=
\{&(1-z)^{\omega_{1}}J^{\varrho_{1}}\}
\{(1-z)^{\omega_{2}-\omega_{1}-\varrho_{1}}J^{\varrho_{2}}\}\dots\\
\{&(1-z)^{\omega_{m-1}-\omega_{m-2}-\varrho_{m-2}}J^{\varrho_{m-1}}\}
\frac{z^{\varrho_m -1}(1-z)^{\omega_m-\omega_{m-1}-\varrho_{m-1}}}{\Gamma(\varrho_m)}
\end{align*}
where $J$ stands for the operation
\[J=\int_0\dots dz\]
This multiple integral can easily be used in the following form:
\[
R\left(z\,\,\vline\begin{tabular}{ccc}$\omega_1$ & $\dots$ & $\omega_m$ \\$\varrho_1$ & $\dots$ & $\varrho_m$ \end{tabular} \right)
=\int_0^z dt_1 \int_0^{t_{1}}dt_2 \dots \int_0^{t_{m-2}}dt_{m-1} \textbf{R}(z|t_1 t_2\dots t_{m-1}) 
\]
\begin{align*}
\textbf{R}(z|t_1 t_2\dots t_{m-1})=&(z-t_1)^{\varrho_1-1}(t_1-t_2)^{\varrho_2-1}\dots(t_{m-2}-t_{m-1})^{\varrho_{m-1}-1}(t_{m-1}^{\varrho_m-1})\times\\
&(1-z)^{\omega_1}(1-t_1)^{\omega_2-\omega_1-\varrho_1}\dots(1-t_{m-2})^{\omega_{m-1}-\omega_{m-2}-\varrho_{m-2}}(1-t_{m-1})^{\omega_{m}-\omega_{m-1}-\varrho_{m-1}}).
\end{align*}
Let's also give a simple Cauchy integral for $R\left( z\,\,\vline\begin{tabular}{ccc}$\omega_1$ & $\dots$ & $\omega_m$ \\$\varrho_1$ & $\dots$ & $\varrho_m$ \end{tabular} \right)$.  It is
\[
R\left(z\,\,\vline \begin{tabular}{ccc}$\omega_1$ & $\dots$ & $\omega_m$ \\ $\varrho_1$ & $\dots$ & $\varrho_m$ \end{tabular} \right)
=
\frac{(-1)^{\sigma-1} \Gamma(\varrho_1) \ldots \Gamma(\varrho_m)}{2 \pi i}
\int_C\frac{(1-z)^\mathfrak{z}d\mathfrak{z}}{\prod^{m}_{k=1}\prod_{h=0}^{\varrho_k-1}(\mathfrak{z}-\omega_k-h)}
\]
which is integrated in the positive direction on a big enough circle about the origin.  Because there is an expansion in decreasing powers
\[
\prod_{k=1}^{m}\prod_{h=0}^{\varrho_k-1}\left(1-\frac{\omega_k+h}{\mathfrak{z}}\right)^{-1}=\sum^\infty_{l=0}b_l\mathfrak{z}^{-l}
\qquad
(b_0=1)
\]
therefore by the theorem of residues we get
\begin{align*}
R\left( z\,\,\vline\begin{tabular}{ccc}$\omega_1$ & $\dots$ & $\omega_m$ \\$\varrho_1$ & $\dots$ & $\varrho_m$ \end{tabular} \right)
&=
(-1)^{\sigma-1}\Gamma(\varrho_1)\dots\Gamma(\varrho_m)
\sum^\infty_{l=0}b_l\frac{(\log(1-z))^{\sigma+l-1}}{\Gamma(\sigma+l)} \\
&=
\frac{\Gamma(\varrho_1)\dots\Gamma(\varrho_m)}{\Gamma(\sigma)}z^{\sigma-1}+\dots
\end{align*}
On the other hand, summing over the residues of the poles of the integrals we get
\begin{align*}
R\left( z\,\,\vline\begin{tabular}{ccc}$\omega_1$ & $\dots$ & $\omega_m$ \\$\varrho_1$ & $\dots$ & $\varrho_m$ \end{tabular} \right)&=\sum^m_{k=1}A_k\left( z\,\,\vline\begin{tabular}{ccc}$\omega_1$ & $\dots$ & $\omega_m$ \\$\varrho_1$ & $\dots$ & $\varrho_m$ \end{tabular} \right)(1-z)^{\omega_k}\\
A_k\left( z\,\,\vline\begin{tabular}{ccc}$\omega_1$ & $\dots$ & $\omega_m$ \\$\varrho_1$ & $\dots$ & $\varrho_m$ \end{tabular} \right)&=(-1)^{\sigma-1}\Gamma(\varrho_1)\dots\Gamma(\varrho_m)\sum^{\varrho_k-1}_{h=0}\frac{(1-z)^h}{\Phi'(\omega_k+h)}\\
\Phi(\mathfrak{z})&=\prod_{k=1}^{m}\prod_{h=0}^{\varrho_k-1}(\mathfrak{z}-\omega_k-h),
\end{align*}
where the polynomial $A_k\left(z\,\,\vline\begin{tabular}{ccc}$\omega_1$ & $\dots$ & $\omega_m$ \\$\varrho_1$ & $\dots$ & $\varrho_m$ \end{tabular} \right)$ is exactly degree $\varrho_k-1$.  Thereby the claims of the definition of $R\left(z\,\,\vline\begin{tabular}{ccc}$\omega_1$ & $\dots$ & $\omega_m$ \\$\varrho_1$ & $\dots$ & $\varrho_m$ \end{tabular} \right)$ are fulfilled.

\subsection*{2.} Using the abbreviation 
\[
F\left(\mathfrak{z}\,\,\vline \begin{tabular}{c}$\omega$ \\$\varrho$ \end{tabular}\right)
=
\prod_{h=0}^{\varrho-1}(\mathfrak{z}-\omega-h)
=\frac{\Gamma(\mathfrak{z}-\omega+1)}{\Gamma(\mathfrak{z}-\omega-\varrho+1)},
\]
we have
\begin{align*}
\Phi(\mathfrak{z})
&=
\prod^m_{k=1}F\left(\mathfrak{z} \,\,\vline \begin{tabular}{c}$\omega_k$ \\$\varrho_k$ \end{tabular} \right),
\\
\Phi'(\omega_k+h)
&=
F'\left( \omega_k+h \,\,\vline \begin{tabular}{c} $\omega_k$ \\ $\varrho_k$ \end{tabular} \right)
\prod^m_{\substack{x=1 \\ x\neq k}}
F\left( \omega_k+h \,\,\vline \begin{tabular}{c} $\omega_x$ \\ $\varrho_x$ \end{tabular} \right)
\end{align*}
for $k=1,2,\dots,m$ and $h=0,1,\dots,\varrho_k-1$.
Further, we have
\[
\frac{\Gamma(\varrho_k)}{F'\left( \omega_k+h\,\,\vline \begin{tabular}{c}$\omega_k$ \\$\varrho_k$ \end{tabular} \right)}=
(-1)^{\varrho_k-h-1}
\begin{pmatrix} \varrho_k-1 \\ h \end{pmatrix},
\]
whereas from the well-known gamma formula we have
\[
\int_G t^x(1+t)^{y-1}dt=\frac{2}{i}\sin\pi x\frac{\Gamma(1+x)\Gamma(y)}{\Gamma(1+x+y)}
\]
for $\Re(y)>0$ where $G$ is the unit-circle in the positive direction and the integral is the principle value.  Thus it follows that for $h=0,1,\dots,\varrho_k-1,\,\,x\neq k$
\[
\frac{\Gamma(\varrho_x)}{F\left( \omega_k+h \,\,\vline \begin{tabular}{c}$\omega_x$ \\$\varrho_x$ \end{tabular} \right)}
=\frac{i(-1)^{\varrho_x-h}}{2\sin(\omega_k-\omega_x)\pi}
\int_G t^{\omega_k-\omega_x+h-\varrho_x}(1+t)^{\varrho_x-1}dt
\]
Thus the $m-1$ variables $t_1,t_2,\dots,t_{k-1},t_{k+1},\dots,t_m$ (integrated respectively in the positive direction on the unit-circles $G_1,G_2,\dots,G_{k-1},G_{k+1},\dots,G_m$ in their planes) are written in the abbreviated form as follows
\[
\int_{G_{1}}dt_{1}\dots\int_{G_{k-1}}dt_{k-1}\int_{G_{k+1}}dt_{k+1}\dots\int_{G_{m}}dt_{m}
=
\int_{(G)}dt.
\]
So now with $Q_k$, the finite and non-zero constant 
\[
Q_k=\prod^m_{\substack{x=1 \\ x\neq k}}\frac{1}{2i\sin(\omega_k-\omega_x)\pi}
\]
we arrive at the following integral formula by means of a simple calculation
\hspace*{-4cm}
\[A_k\left(z\,\,\vline\begin{tabular}{ccc}$\omega_1$ & $\dots$ & $\omega_m$ \\$\varrho_1$ & $\dots$ & $\varrho_m$ \end{tabular} \right)
=Q_k(-1)^{m-1}\int_{(G)}\left(1+(-1)^m t_1\dots t_{k-1}t_{k+1}\dots t_m(1-z)\right)^{\varrho_k-1}
\prod^m_{\substack{x=1 \\ x\neq k}}t_x^{\omega_k-\omega_x-\varrho_x}(1+t_x)^{\varrho_x-1}
dt
.
\]
\subsection*{3.} We define the symbol $\delta_{hk}$ for $(h,k=1,2,\dots,m)$ as follows
\[
\delta_{hk}=\left \{ \begin{tabular}{ccc} $1$&for&$h=k$ \\ $0$&for&$h\neq k$ \end{tabular}\right.
\]
and for $(h,k=1,2,\dots,m)$ set 
\begin{align*}
R_h\left(z\,\,\vline\begin{tabular}{ccc}$\omega_1$ & $\dots$ & $\omega_m$ \\$\varrho_1$ & $\dots$ & $\varrho_m$ \end{tabular} \right)
&=
R\left(z\,\,\vline\begin{tabular}{ccc}$\omega_1$ & $\dots$ & $\omega_m$ \\$\varrho_1+\delta_{h1}$ & $\dots$ & $\varrho_m+\delta_{hm}$ \end{tabular} \right)
\\
A_{hk}\left(z\,\,\vline\begin{tabular}{ccc}$\omega_1$ & $\dots$ & $\omega_m$ \\$\varrho_1$ & $\dots$ & $\varrho_m$ \end{tabular} \right)
&=
A_{k}\left(z\,\,\vline\begin{tabular}{ccc}$\omega_1$ & $\dots$ & $\omega_m$ \\$\varrho_1+\delta_{h1}$ & $\dots$ & $\varrho_m+\delta_{hm}$ \end{tabular} \right)
.\end{align*}

Thus between the determinant 
\[
\Delta\left(z\,\,\vline\begin{tabular}{ccc}$\omega_1$ & $\dots$ & $\omega_m$ \\$\varrho_1$ & $\dots$ & $\varrho_m$ \end{tabular} \right)
=
\left|
A_{hk}\left(z\,\,\vline\begin{tabular}{ccc}$\omega_1$ & $\dots$ & $\omega_m$ \\$\varrho_1$ & $\dots$ & $\varrho_m$ \end{tabular} \right)
\right|
\]
and the minor
\[
\Delta_{hk}\left(z\,\,\vline\begin{tabular}{ccc}$\omega_1$ & $\dots$ & $\omega_m$ \\$\varrho_1$ & $\dots$ & $\varrho_m$ \end{tabular} \right)
=
\left|
A_{h'k'}\left(z\,\,\vline\begin{tabular}{ccc}$\omega_1$ & $\dots$ & $\omega_m$ \\$\varrho_1$ & $\dots$ & $\varrho_m$ \end{tabular} \right)
\right|_{\substack{h'\neq h \\ k'\neq k}}
\]
there is the identity
\[
\Delta\left(z\,\,\vline\begin{tabular}{ccc}$\omega_1$ & $\dots$ & $\omega_m$ \\$\varrho_1$ & $\dots$ & $\varrho_m$ \end{tabular} \right)(1-z)^{\omega_k}
=
\sum_{h=1}^m(-1)^{h+k}
\Delta_{hk}\left(z\,\,\vline\begin{tabular}{ccc}$\omega_1$ & $\dots$ & $\omega_m$ \\$\varrho_1$ & $\dots$ & $\varrho_m$ \end{tabular} \right)
R_h\left(z\,\,\vline\begin{tabular}{ccc}$\omega_1$ & $\dots$ & $\omega_m$ \\$\varrho_1$ & $\dots$ & $\varrho_m$ \end{tabular} \right)
\]
for $(k=1,2,\dots,m)$.
Thus $\Delta\left(z\,\,\vline\begin{tabular}{ccc}$\omega_1$ & $\dots$ & $\omega_m$ \\$\varrho_1$ & $\dots$ & $\varrho_m$ \end{tabular} \right)$ has a root of order $\sigma$ at $z=0$; Now since it is cleary also a polynomial of order exactly $\sigma$, the following must hold
\[
\Delta\left(z\,\,\vline\begin{tabular}{ccc}$\omega_1$ & $\dots$ & $\omega_m$ \\$\varrho_1$ & $\dots$ & $\varrho_m$ \end{tabular} \right)
=
\delta\left(\begin{tabular}{ccc}$\omega_1$ & $\dots$ & $\omega_m$ \\$\varrho_1$ & $\dots$ & $\varrho_m$ \end{tabular} \right)z^\sigma,
\]
in which the constant $\delta\left(\begin{tabular}{ccc}$\omega_1$ & $\dots$ & $\omega_m$ \\$\varrho_1$ & $\dots$ & $\varrho_m$ \end{tabular} \right)$ is independent of $z$; Furthmore since $A_{hk}\left(z\,\,\vline\begin{tabular}{ccc}$\omega_1$ & $\dots$ & $\omega_m$ \\$\varrho_1$ & $\dots$ & $\varrho_m$ \end{tabular} \right)$ is a polynomial in $z$ of degree exactly $\varrho_k+\delta_{hk}-1$, it follows that $\delta\left(\begin{tabular}{ccc}$\omega_1$ & $\dots$ & $\omega_m$ \\$\varrho_1$ & $\dots$ & $\varrho_m$ \end{tabular} \right)$ is not zero; 
\textit{Thus the determinant}
\[
\Delta\left(z\,\,\vline\begin{tabular}{ccc}$\omega_1$ & $\dots$ & $\omega_m$ \\$\varrho_1$ & $\dots$ & $\varrho_m$ \end{tabular} \right)
\]
\textit{vanishes if and only if $z=0$\footnote{Carrying out the calculation yields the value
\[
\delta\left(\begin{tabular}{ccc}$\omega_1$ & $\dots$ & $\omega_m$ \\$\varrho_1$ & $\dots$ & $\varrho_m$ \end{tabular} \right)
=\mp\prod^m_{\substack{h,k=1\\h\neq k}}\frac{\Gamma(\omega_h-\omega_k)\Gamma(\varrho_k)}{\Gamma(\varrho_k+\omega_h-\omega_k)}\neq 0
\]
see (8).
}
.}

\section*{II.}
\subsection*{4.} Let $n$ be a natural number such that $n\geq3$ and $n\geq m\geq2$ and
\[
R_h(z)=R_h\left(z\,\,\vline\begin{tabular}{ccc}$0\,\frac{1}{n}$&$\dots$&$\frac{m-1}{n}$ \\$\varrho\,\varrho$&$\dots$&$\varrho$ \end{tabular} \right)
,\,\,\,\,\,\,\,
A_{hk}(z)=A_{hk}\left(z\,\,\vline\begin{tabular}{ccc}$0\,\frac{1}{n}$&$\dots$&$\frac{m-1}{n}$ \\$\varrho\,\varrho$&$\dots$&$\varrho$ \end{tabular} \right)
\] 
for $h,k=1,2,\dots,m$, so that 
\[
R_h(z)=\sum^m_{k=1}A_{hk}(z)(1-z)^{\frac{k-1}{n}}
\]
and so that $R_h(z)$ has root of order $m\varrho$ at $z=0$.  With the new variables
\[
x=x(z)=(1-z)^{\frac{1}{n}},\,\,\,\,\,\,x(0)=1
\]
we compose and rewrite the previous functions in the following manner
\[
\mathfrak{R}_h(x)=R_h(1-x^n),\,\,\,\,\,\,\mathfrak{A}_{hk}(x^n)=A_{hk}(1-x^n)
.
\]
The neighborhood of $z=0$ is mapped to the neighborhood of $x=1$; $\mathfrak{R}_h(x)$ has therefore at $x=1$ a root of order $m\varrho$.  Setting
\[
\mathfrak{S}_h(x)=(x-1)^{-m\varrho}\mathfrak{R}_h(x),
\]
we then have that $\mathfrak{S}_h(x)$ is regular in a neighborhood of $x=1$; Thus we can easily see that $\mathfrak{S}_h(x)$ is a polynomial.

Introducing yet another independent variable $y$ and setting
\begin{align*}
\mathfrak{T}_h(x\,\,y)&=\sum^m_{k=1}\mathfrak{A}_{hk}(x^n)
\frac{y^{k-1}-x^{k-1}}{y-x},
\\
\mathfrak{U}_h(x\,\,y)&=\sum^m_{k=1}\mathfrak{A}_{hk}(x^n)y^{k-1},
\end{align*}
 we then have the identity  
\[
\mathfrak{U}_h(x\,\,y)=(x-1)^{m\varrho}\mathfrak{S}_h(x)+(y-x)\mathfrak{T}_h(x\,\,y)
\]
for $(h=1,2,\dots,m)$.  From subsection 3. The determinant
\[\left|\mathfrak{A}_{hk}(x^n)\right|\]
is  non-zero, if $x$ is not an $n$-th root of unity.  So following from this condition for every value of $y$ at least one of the $m$ numbers
\[
\mathfrak{U}_h(x,y)
\]
for $(h=1,2,\dots,m)$ is non-zero.

\subsection*{5.} Let $a,b$ be two natural numbers such that
\[\xi =\sqrt[n]{\frac{a}{b}}\]
is an algebraic number of degree exactly $n$. Any two rational numbers with positive denominators $\frac{p_1}{q_1}$ and $\frac{p_2}{q_2}$, that satisfy the inequalities
$
\left(\frac{2}{3}\right)^{\frac{1}{n}}\xi\leq\frac{p_1}{q_1}\leq\left(\frac{3}{2}\right)^{\frac{1}{n}}\xi
,\,\,\,\,\,\,
\left(\frac{2}{3}\right)^{\frac{1}{n}}\xi\leq\frac{p_2}{q_2}\leq\left(\frac{3}{2}\right)^{\frac{1}{n}}\xi
$
may be used to assign the following values to $x$ and $y$
\[x=\frac{q_1\xi}{p_1},\,\,\,\,\,\,y=\frac{q_1p_2}{p_1q_2}.\]
As long as $x$ is not an $n$-th root of unity then at least one of the $m$ $\mathfrak{U}_h(x,y)$'s is non-zero.  So for some $h_0$ we have
\[\mathfrak{U}_{h_0}(x,y)\neq0.\]
Clearly $\mathfrak{U}_{h_0}(x,y)$ is a rational number whose denominator can be estimated to its upper limits.

It was claimed earlier that
\[
\varrho(\varrho-1)^{m-1}A_{hk}(z)=(-1)^{m\varrho}\sum^{\varrho+\delta_{hk}-1}_{l=0}
\frac{(\varrho+1)!(\varrho!)^{m-1}}{\Phi^{'}_{h}\left(\frac{k-1}{n}+l\right)}(1-z)^l
\]
with
\[
\frac{(\varrho+1)!(\varrho!)^{m-1}}{\Phi^{'}_{h}\left(\frac{k-1}{n}+l\right)}
=\frac{(\varrho+\delta_{hk})!}{F^{'}\left(\frac{k-1}{n}+l\,\,\vline\,\,\substack{\frac{k-1}{n}\\ \varrho+\delta_{hk}} \right)}
\prod^m_{\substack{x=1\\x\neq k}}
\frac{(\varrho+\delta_{hx})!}{F^{'}\left(\frac{k-1}{n}+l\,\,\vline\,\,\substack{\frac{x-1}{n}\\ \varrho+\delta_{hx}} \right)},
\]
in which
\[
\frac{(\varrho+\delta_{hk})!}{F^{'}\left(\frac{k-1}{n}+l\,\,\vline\,\,\substack{\frac{k-1}{n}\\ \varrho+\delta_{hk}} \right)}
\]
ends up being entirely rational, where on the other hand for $x\neq k$ we have
\[
\frac{(\varrho+\delta_{hx})!}{F\left(\frac{k-1}{n}+l\,\,\vline\,\,\substack{\frac{x-1}{n}\\ \varrho+\delta_{hx}} \right)}
=(-1)^{\varrho+\delta_{hx}}
\frac{(n)(2n)\dots (( \varrho+\delta_{hx})n)}{(n+K)(2n+K)\dots((\varrho+\delta_{hx})n+K)}
\]
with $K=-nl+x-k-n$.

According to a theorem of Maier the lowest common denominator of the coefficients of all the polynomials $A_{hk}(z)$ must be smaller than the $\varrho$-th power of a constant that depends only on $n$ and $m$\footnote{See the works (9) of Maier and (7) of Siegel, where the proofs are carried out}.  On account of
\[
\mathfrak{U}_{h_0}(xy)=\sum^m_{k=1}A_{h_0k}\left(1-\frac{aq^n_1}{bp^n_1}\right)
\left(\frac{q_1p_2}{p_1q_2}\right)^{k-1}
\]
we have the denominator of the rational number $\mathfrak{U}_{h_0}(xy)$.  Therefore through multiplication with $b^\varrho p_1^{n\varrho+m-1}q_2^{m-1}$, and since $\mathfrak{U}_{h_0}(xy)$ is not equal to zero, there exists the inequality
\[
\left|\mathfrak{U}_{h_0}(xy)\right|^{-1}\leq c_1^\varrho p_1^{n\varrho+m-1}q_2^{m-1} 
\]
with positive constant $c_1$ that depends only on $n$,$m$ and $b$.

Further from subsection 1. we have
\[
R_h(z)=\int_0^{z}dt_{1}\int_0^{t_1}dt_{2}\int_0^{t_{m-2}}dt_{m-1}
\frac{(z-t_1)^{\varrho+\delta_{h1}-1}(t_1-t_2)^{\varrho+\delta_{h2}-1}\dots(t_{m-2}-t_{m-1})^{\varrho+\delta_{h(m-1)}-1}(t_{m-1})^{\varrho+\delta_{hm}}}
{(1-t_{1})^{\varrho+\delta_{h1}-\frac{1}{n}}(1-t_{2})^{\varrho+\delta_{h2}-\frac{1}{n}}\dots(1-t_{m-1})^{\varrho+\delta_{h(m-1)}-\frac{1}{n}}}
\]
and with the new variables of integration $t_k=zu_k$ $(k=1,2,\dots,m)$
\begin{align*}
R_h(z)=&z^{m\varrho}\int_0^{1}du_{1}\int_0^{u_1}du_{2}\int_0^{u_{m-2}}du_{m-1}
\frac{(1-u_1)^{\varrho+\delta_{h1}-1}(u_1-u_2)^{\varrho+\delta_{h2}-1}\dots(u_{m-2}-u_{m-1})^{\varrho+\delta_{h(m-1)}-1}(u_{m-1})^{\varrho+\delta_{hm}}}
{(1-zu_{1})^{\varrho+\delta_{h1}-\frac{1}{n}}(1-zu_{2})^{\varrho+\delta_{h2}-\frac{1}{n}}\dots(1-zu_{m-1})^{\varrho+\delta_{h(m-1)}-\frac{1}{n}}}
\\
=&z^{m\varrho}J
\end{align*}
whereby due to 
\[
z=1-x^n\leq\frac{1}{2}
\]
the factors of the denominators are greater than $\frac{1}{2}$.  Further, since
\[
\mathfrak{S}_h(x)=\left(-\frac{1-x^n}{1-x}\right)^{m\varrho}J,\,\,\,\,\,\,
\left|\frac{1-x^n}{1-x}\right|=\left|1+x+\dots+x^{n-1}\right|\leq\frac{3n}{2},
\]
we have the following inequality
\[\left|\mathfrak{S}_h(x)\right|\leq c^{\varrho}_2,\]
wherein the positive constant $c_2$ depends only on $n$ and $m$.

Finally, it follows from the integral formula in subsection 2. that
\[
(-1)^{m-1}A_{hk}(z)=\prod^m_{\substack{x=1\\x\neq k}}\frac{1}{2i\sin\frac{k-x}{n}\pi}
\int_{(G)}\Big(1+(-1)^m(1-z)\prod^m_{\substack{x=1\\x\neq k}}t_x\Big)^{\varrho+\delta_{hk}-1}\prod^m_{\substack{x=1\\x\neq k}}t_x^{\frac{k-x}{n}-\varrho-\delta_{hx}}(1+t_x)^{\varrho+\delta_{hk}-1}dt
\]
and from the definition of $\mathfrak{T}_h(xy)$ that
\[\left|\mathfrak{T}_h(xy)\right|\leq c^{\varrho}_3\]
wherein the positive constant $c_3$ depends only on $n$ and $m$.

All members of the identity
\[
\mathfrak{U}_{h_0}(xy)=\left(\frac{q_1}{p_1}\right)^{m\varrho}\mathfrak{S}_{h_0}(x)\left(\xi-\frac{p_1}{q_1}\right)^{m\varrho}+\frac{q_1}{p_1}\mathfrak{T}_{h_0}(xy)\left(\frac{p_2}{q_2}-\xi\right)
\]
have their values derived either above or below and from them follows the existence of two positive constants $c_4$ and $c_5$, which in turn only depend on $n$,$m$ and $\xi$ so that the sum of the two numbers
\[\vartheta_1=c^{\varrho}_4q_1^{n\varrho+m-1}q_2^{m-1}\left|\xi-\frac{p_1}{q_1}\right|^{m\varrho}
,\,\,\,\,\,\,
\vartheta_2=c^{\varrho}_5q_1^{n\varrho+m-1}q_2^{m-1}\left|\xi-\frac{p_2}{q_2}\right|
\]
is greater than two and at least one of the numbers is also greater than one.

\subsection*{6.} Now we easily succeed in proving the Thue-Siegel theorem for the specific algebraic numbers $\xi$ \footnote{Refer to works (1) though (7)}:\\``With $\varepsilon$ an arbitrary constant and $m$ an arbitrary natural number it follows that the inequality 
\[
\left|\xi -\frac{p}{q}\right|\leq q^{-(\frac{n}{m}+m-1)-\epsilon}
\]
has only finitely many rational solutions $\frac{p}{q}$ with positive denominator.''

It suffices to limit the proof of $m$ to the numbers $2,3,\dots,n$.  Only those solutions of the previous inequality need be considered that also satisfy the following inequality
\[
\left(\frac{2}{3}\right)^{\frac{1}{n}}\xi\leq\frac{p}{q}\leq\left(\frac{3}{2}\right)^{\frac{1}{n}}\xi
;
\]
which is a consequence of the first if the denominator $q$ is sufficiently large.  For any two such rational numbers $\frac{p_1}{q_1}$ and $\frac{p_2}{q_2}$, the natural number $\varrho$ is determined by the condition
\[
q_1^{m(\varrho-1)}<q_2\leq q_1^{m\varrho}.
\]
A simple calculation yields the inequalities
\[
\vartheta_1\leq q_1^{m-1} q_2^{-\frac{\varepsilon}{2}}\left(c_4q_1^{-\frac{m\varepsilon}{2}}\right)^\varrho;
\,\,\,\,\,\,
\vartheta_2\leq q_1^{n+m+m\varepsilon-1} q_2^{-\frac{\varepsilon}{2}}\left(c_5q_1^{-\frac{m\varepsilon}{2}}\right)^\varrho
.
\]
Now if there were infinitely many solutions for the inequality
\[
\left|\xi -\frac{p}{q}\right|\leq q^{-(\frac{n}{m}+m-1)-\epsilon}
\]
then it could be show that
\begin{align*}
&\left(\frac{2}{3}\right)^{\frac{1}{n}}\xi\leq\frac{p_1}{q_1}\leq\left(\frac{3}{2}\right)^{\frac{1}{n}}\xi,\,\,\,\,\,\,q_1\geq\max\left(c_4^{\frac{2}{m\varepsilon}},c_5^{\frac{2}{m\varepsilon}}\right),
\\
&\left(\frac{2}{3}\right)^{\frac{1}{n}}\xi\leq\frac{p_2}{q_2}\leq\left(\frac{3}{2}\right)^{\frac{1}{n}}\xi,\,\,\,\,\,\,q_2\geq\max\left(q_1^{\frac{2}{\varepsilon}(m-1)},q_1^{\frac{2}{\varepsilon}(n+m+m\varepsilon-1)}\right);
\end{align*}
but this in turn would mean that $\vartheta_1\leq1$ and $\vartheta_2\leq1$, which contradicts what has been shown above.

\vspace{5mm}
G{\"o}ttingen, February 4\textsuperscript{th} 1931.
\vspace{5mm}
\section*{Bibliography}
\begin{tabular}{p{.16cm} p{1.65cm} p{14cm}}
(1) & A. Thue:&Bemerkung {\"u}ber gewisse N{\"a}hrungsbr{\"u}che algebraischer Zahlen, {\em Videnskapsselskapets-Skrifter Christiania} (1908).\\
(2) & &Om en general i store hele tal ul{\o}sbar ligning, {\em Videnskapsselskapets-Skrifter Christiania}  (1908).\\
(3) & &{\"U}ber Ann{\"a}herungswerte algebraischer Zahlen, {\em Journal f{\"u}r die reine und angewandte Mathematik} \textbf{135} (1909).\\
(4) & C. Siegel:&Approximation algebraischer Zahlen, {\em Math. Zeitschr.} \textbf{10} (1921).\\
(5) & &N{\"a}herungswerte algebraischer Zahlen, {\em Math. Annalen} \textbf{84} (1921).\\
(6) & &{\"U}ber den Thueschen Satz, {\em Videnskapsselskapets-Skrifter Christiania} (1922).\\
(7) & &{\"U}ber einige Anwendungen diophantischer Approximationen, {\em Abh. d. Preu{\ss}. Akad. d. Wissensch.} (1929).\\
(8) & K. Mahler:& Zur Approximation der Exponentialfunktion und des Logarithmus, Journal f{\"u}r die reine und angewandte Mathematik. (At the printer).\\
(9) & W. Maier:&Potenzreihen irrationalen Grenzwertes, {\em Journal f{\"u}r die reine und angewandte Mathematik} \textbf{156} (1925).\\
\end{tabular}
\vspace{.5mm}

Further work of Thue on Diophantine approximation is cited in the bibliography of (4).
\vspace{7mm}
\center{(Submitted on February 16\textsuperscript{th} 1931.)}



\end{document}